\documentclass[11pt,twoside,A4paper]{article} 
\usepackage{amssymb,amsthm,amsmath,amscd,amsfonts,txfonts}
\usepackage{xcolor}

\usepackage[tmargin=1.5cm,bmargin=1.5cm,lmargin=1.7cm,rmargin=1.7cm]{geometry}

\theoremstyle{plain}

\newtheorem{proposition}{Proposition}[section]
\newtheorem{theorem}{Theorem}[section]
\newtheorem{corollary}{Corollary}[section]

\theoremstyle{definition}
\newtheorem{definition}{Definition}[section]

\theoremstyle{remark}
\newtheorem{remark}{Remark}[section]

\title{\bfseries\scshape{Derivations, Centroids and Iner-Derivations of $\textbf{5}$-Dimensional Nilpotent Complex Associative Algebras}}
\author{
\bfseries\scshape Ahmed Zahari ABDOU \thanks{e-mail address: zaharymaths@gmail.com}\\
Universit\'{e} de Haute Alsace,\\
 IRIMAS-D\'{e}partement de Math\'{e}matiques.}

\date{} 

\begin{document} 
\maketitle

\tableofcontents

\begin{abstract} 
\noindent This study focuses on the analysis of derivations, centroids, and inner derivations of 5-dimensional complex nilpotent associative algebras. It presents the classification of these algebras of dimension less than five, as well as the classifications of their corresponding derivations, centroids and iner-derivations.
\end{abstract} 

\noindent
{\bf Mathematics Subject Classification 2020:} 17A40, 17A99\\
\noindent
{\bf Keywords:} Associative algebra, classification, derivation, centroid, inner-derivative. 

\section{Introduction}

The study of finite-dimensional associative algebras, and nilpotent algebras in particular, plays a key role in pure and applied mathematics. These structures appear naturally in several fields, such as deformation theory \cite{G1}, homological algebra, quantum mechanics and noncommutative geometry \cite{N}. Associative nilpotent algebras over complex numbers are distinguished by their complex structure and the various behaviors they exhibit in their classification, particularly as their dimension increases. Motivated by their importance in the analysis of rigidity properties, deformations and isomorphism classes, this study focuses on complex nilpotent associative algebras of dimension five, whose classification reveals a great structural richness.

A five-dimensional complex nilpotent associative algebra is a five-dimensional complex vector space equipped with an associative bilinear multiplication, where repeated products of elements eventually vanish. Formally, an associative algebra \( A \) is called nilpotent if there exists an integer \( k \) such that any product of \( k \) elements (in any associative arrangement) is zero, i.e., \( A^k = 0 \), where \( A^k \) denotes the ideal generated by all products of \( k \) elements. In the context of associative algebras, nilpotency implies a highly constrained multiplication structure, rich in degeneracies, yet capable of revealing deep mathematical properties \cite{M}.

When analyzing the internal and external symmetries of these algebras, three key concepts are essential: derivations, centroids, and inner derivations. Let \( (A, \cdot) \) be an associative algebra. A linear map \( d : A \to A \) is called a derivation if it satisfies the Leibniz rule:
\[
d(x \cdot y) = d(x) \cdot y + x \cdot d(y)
\]
for all \( x, y \in A \). Centroids also play a significant role in understanding central extensions and module structures over the algebra. An algebra with a large centroid may admit numerous non-trivial deformations and central extensions, thus exhibiting a more flexible algebraic behavior \cite{J}.

In associative algebras, an inner derivation associated with an element \( a \in A \) is defined via the commutator:
\[
\delta_a(x) = [a, x] = a x - x a.
\]
Examining the space of inner derivations relative to the full derivation algebra \( \text{Der}(A) \) provides insights into the rigidity or flexibility of the algebra’s internal symmetries. In many cases, especially for nilpotent algebras of low dimension, \( \text{Der}(A) \) strictly contains the space of inner derivations, thus revealing the presence of external derivations \cite{V}.

Dimension five marks a crucial threshold in the study of nilpotent associative algebras. While classification results for algebras of dimension 2 to 4 are relatively well-established and manageable, the case of dimension five is significantly more complex, revealing a broader variety of non-isomorphic algebraic structures even within the class of nilpotent algebras \cite{GF}.

The structure of this paper is as follows:
\begin{itemize}
	\item \textbf{Section 2}: Provides basic definitions and properties of associative algebras, derivations, centroids, and inner derivations;
	\item \textbf{Section 3}: Introduces a step-by-step method for computing derivations of five-dimensional associative algebras;
	\item \textbf{Section 4}: Describes a method for computing centroids of five-dimensional associative algebras;
	\item \textbf{Section 5}: Presents a method for calculating inner derivations of five-dimensional associative algebras.
\end{itemize}

This work focuses on a detailed analysis of derivations, centroids, and inner derivations in five-dimensional complex nilpotent associative algebras. It begins by reviewing the classification of nilpotent associative algebras of dimension less than five, which provides a foundation for understanding the structures in dimension five \cite{M}. A systematic exploration follows of the derivations, centroids, and inner derivations associated with each type of classified algebra.

The study of these structures is driven by their essential role in the investigation of algebra deformations \cite{G1}, as well as their significance in understanding rigidity properties and isomorphism classes of nilpotent algebras \cite{V}. In particular, centroids are central to algebra classification and the analysis of inner automorphisms \cite{J, s, Ab}.

It is a continuation of recent research on the classification of nilpotent algebras, with particular emphasis on 2-nilpotent algebras. We were inspired by the geometric approach presented in \cite{MIY, u}, which provides a rigorous framework for the study of varieties of nilpotent algebras using algebraic geometry and the action of algebraic groups. Building on this methodology, our main goal is to enhance the understanding of these structures by emphasizing their internal properties, particularly their fundamental invariants, while exploring connections between associative algebra theory, cohomology, and linear algebra. This approach offers a fresh perspective on the fine structure of these algebraic objects, enriching the existing classification, particularly in small dimensions.

\section{Classification}
\begin{definition}\label{d1}
Let $\textbf{A}$ be a vector space over a field $\mathbb{K}$.

\begin{itemize}
  \item[(i)] An associative algebra structure on $\textbf{A}$ is a bilinear map
  \[
  \cdot : \textbf{A} \otimes \textbf{A} \rightarrow \textbf{A}, \quad (x, y) \mapsto x \cdot y,
  \]
  such that the multiplication is associative, that is,
  \[
  (x \cdot y) \cdot z = x \cdot (y \cdot z), \quad \forall x, y, z \in \textbf{A}.
  \]

  \item[(ii)] The center of the associative algebra $\textbf{A}$ is defined as:
  \[
  Z_{\textbf{A}} = \left\{ x \in \textbf{A} \,\middle|\, x \cdot y = y \cdot x, \ \forall y \in \textbf{A} \right\}.
  \]
\end{itemize}
\end{definition}

\begin{definition}\label{d2}
 Let $(\textbf{A}, \cdot)$ be a associative algebra. A linear map $d : \textbf{A}\rightarrow\textbf{A}$ is said to be 
\begin{itemize}
  \item[(i)] a {\bf derivation}  if \hspace{0,25cm} 
$d(x\cdot y) = d(x)\cdot y + x\cdot d(y).$
  \item[(ii)] an {\bf element of centroid} \hspace{0,25cm}  if
$\phi(x\cdot y)=\phi(x)\cdot y=x\cdot\phi(y)$,
\end{itemize}
for any $x, y\in\textbf{A}$.
\end{definition}

By utilizing the classification results \cite{MIY} and solving the equations, the derivations and automorphism groups can be determined and summarized.

\begin{theorem}
The isomorphism class of $\textbf{5}$-dimensional complex nilpotent (non-2-step nilpotent) non-commutative associative algebras given by  the following representatives.
\end{theorem}
\small\begin{center}
    \begin{tabular}{||c||c||c||c||c||c||c||}
\hline
$\textbf{A}_i$&Multiplications&Center&$\textbf{A}_i$&Multiplications&Center.
\\ \hline 
$\textbf{A}_1$&
$\begin{array}{ll}  
e_1 \cdot e_1=e_2\\
e_1 \cdot e_2=e_4\\
e_1  \cdot e_3=e_4\\
\end{array}$
$\begin{array}{ll}  
e_2  \cdot e_1=e_4\\
e_3  \cdot e_3=e_4
\end{array}$
&
$\textless e_2,e_4,e_5\textgreater$
&
$\textbf{A}_2$
&
$\begin{array}{ll}
e_1 \cdot e_1=e_2\\
e_1e_2=e_4\\
\end{array}$
$\begin{array}{ll}  
e_1  \cdot e_3=e_4\\
e_2  \cdot e_1=e_4
\end{array}$
&
$\textless e_2,e_4,e_5\textgreater$
\\ \hline 
$\textbf{A}_3$&
$\begin{array}{ll}  
e_1 \cdot e_1=e_2\\
e_1  \cdot e_2=e_3\\
e_1  \cdot e_3=e_4\\
e_1  \cdot e_5=e_4\\
\end{array}$
$\begin{array}{ll}  
e_2 \cdot e_1=e_3\\
e_2  \cdot e_2=e_4\\
e_3 \cdot e_1=e_4
\end{array}$
&
$\textless e_2,e_4,e_5\textgreater$
&
$\textbf{A}_4$
&
$\begin{array}{ll}
e_1  \cdot e_1=e_2\\
e_1  \cdot e_2=e_4\\
e_1  \cdot e_3=e_4\\
e_1 \cdot e_5=e_4\\
\end{array}$
$\begin{array}{ll}  
e_2  \cdot e_1=e_3\\
e_2 \cdot e_2=e_4\\
e_3 \cdot e_1=e_4\\
e_5e_5=e_4
\end{array}$
&
$\textless e_2,e_4,e_5\textgreater$
\\ \hline 
$\textbf{A}_5$&
$\begin{array}{ll}  
e_1 \cdot e_1=e_2\\
e_1 \cdot e_2=e_3\\
e_1 \cdot e_4=e_3
\end{array}$
$\begin{array}{ll}  
e_2 \cdot e_1=e_3\\
e_4 \cdot e_5=e_3\\
e_5 \cdot e_4=e_3
\end{array}$
&
$\textless e_2,e_5\textgreater$
&
$\textbf{A}_6$
&
$\begin{array}{ll}
e_1 \cdot e_1=e_2\\
e_1 \cdot e_2=e_3\\
e_1 \cdot e_4=e_3
\end{array}$
$\begin{array}{ll}  
e_2 \cdot e_1=e_3\\
e_5 \cdot e_5=e_3
\end{array}$
&
$\textless e_2,e_3,e_5\textgreater$
\\ \hline 
$\textbf{A}_7$&
$\begin{array}{ll}  
e_1 \cdot e_1=e_2\\
e_1 \cdot e_2=e_3\\
e_1 \cdot e_4=e_3
\end{array}$
$\begin{array}{ll}  
e_2 \cdot e_1=e_3\\
e_4 \cdot e_4=e_3\\
e_5 \cdot e_5=e_3
\end{array}$
&
$\textless e_2,e_3,e_5\textgreater$
&
$\textbf{A}_8$
&
$\begin{array}{ll}
e_1 \cdot e_1=e_2\\
e_1 \cdot e_2=e_3\\
e_2 \cdot e_1=e_3
\end{array}$
$\begin{array}{ll}  
e_4 \cdot e_5=e_3\\
e_5 \cdot e_4=-e_3\\
e_5 \cdot e_5=e_3
\end{array}$
&
$\textless e_1,e_2,e_3\textgreater$
\\ \hline
$\textbf{A}_9$&
$\begin{array}{ll}  
e_1 \cdot e_1=e_2\\
e_1 \cdot e_2=e_3\\
e_2 \cdot e_1=e_3
\end{array}$
$\begin{array}{ll}  
e_2 \cdot e_1=e_3\\
e_4 \cdot e_5=e_3\\
e_5 \cdot e_4=\alpha e_3
\end{array}$
&
$\textless e_1,e_2,e_3\textgreater$
&
$\textbf{A}_{10}$
&
$\begin{array}{ll}
e_1 \cdot e_1=e_2\\
e_1 \cdot e_2=e_3\\
e_2 \cdot e_1=e_3
\end{array}$
$\begin{array}{ll}  
e_2 \cdot e_1=e_3\\
e_4 \cdot e_1=e_5
\end{array}$
&
$\textless e_2,e_3,e_5\textgreater$
\\ \hline 
$\textbf{A}_{11}$&
$\begin{array}{ll}  
e_1 \cdot e_1=e_2\\
e_1 \cdot e_2=e_3\\
e_2 \cdot e_1=e_3
\end{array}$
$\begin{array}{ll}  
e_4 \cdot e_1=e_3\\
e_4 \cdot e_4= e_3
\end{array}$
&
$\textless e_2,e_3,e_5\textgreater$
&
$\textbf{A}_{12}$
&
$\begin{array}{ll}
e_1 \cdot e_1=e_2\\
e_1 \cdot e_2=e_3\\
e_2 \cdot e_1=e_3\\
\end{array}$
$\begin{array}{ll}  
e_4 \cdot e_1=e_5\\
e_4 \cdot e_2=e_3\\
e_5 \cdot e_1=e_3
\end{array}$
&
$\textless e_3\textgreater$
\\ \hline 
$\textbf{A}_{13}$&
$\begin{array}{ll}  
e_1 \cdot e_1=e_2\\
e_1 \cdot e_2=e_3\\
e_2 \cdot e_1=e_3\\
e_4 \cdot e_1=e_5
\end{array}$
$\begin{array}{ll}  
e_4 \cdot e_1=e_3\\
e_4 \cdot e_4=e_3\\
e_5 \cdot e_1=e_3
\end{array}$
&
$\textless e_3\textgreater$
&
$\textbf{A}_{14}$
&
$\begin{array}{ll}
e_1 \cdot e_1=e_2\\
e_1 \cdot e_2=e_3\\
e_1 \cdot e_4=e_5
\end{array}$
$\begin{array}{ll}  
e_2 \cdot e_1=e_3\\
e_4 \cdot e_1=e_3+e_5
\end{array}$
&
$\textless e_3\textgreater$
\\ \hline 
$\textbf{A}_{15}$&
$\begin{array}{ll}  
e_1 \cdot e_1=e_2\\
e_1 \cdot e_2=e_3\\
e_1 \cdot e_4=e_5
\end{array}$
$\begin{array}{ll}  
e_2 \cdot e_1=e_3\\
e_4 \cdot e_1=e_3+e_5\\
e_4 \cdot e_4=e_3
\end{array}$
&
$\textless e_2,e_3,e_5\textgreater$
&
$\textbf{A}_{16}$
&
$\begin{array}{ll}
e_1 \cdot e_1=e_2\\
e_1 \cdot e_2=e_3\\
e_1 \cdot e_4=e_5
\end{array}$
$\begin{array}{ll}  
e_2 \cdot e_1=e_3\\
e_4 \cdot e_1=\alpha e_5
\end{array}$
&
$\textless e_2,e_3,e_5\textgreater$
\\ \hline 
$\textbf{A}_{17}$&
$\begin{array}{ll}  
e_1 \cdot e_1=e_2\\
e_1 \cdot e_2=e_3\\
e_1 \cdot e_4=e_5
\end{array}$
$\begin{array}{ll}  
e_2 \cdot e_1=e_3\\
e_4 \cdot e_1=\alpha e_5\\
e_4 \cdot e_4=e_3
\end{array}$
&
$\textless e_2,e_3,e_5\textgreater$
&
$\textbf{A}_{18}$
&
$\begin{array}{ll}
e_1 \cdot e_1=e_2\\
e_1 \cdot e_2=e_3\\
e_2 \cdot e_1=e_3
\end{array}$
$\begin{array}{ll}  
e_4 \cdot e_1=e_3\\
e_4 \cdot e_4=e_5
\end{array}$
&
$\textless e_2,e_3,e_5\textgreater$
\\ \hline  
$\textbf{A}_{19}$&
$\begin{array}{ll}  
e_1 \cdot e_1=e_2\\
e_1 \cdot e_2=e_3\\
e_1 \cdot e_4=e_3
\end{array}$
$\begin{array}{ll}  
e_2 \cdot e_1=e_3\\
e_4 \cdot e_4=e_5
\end{array}$
&
$\textless e_2,e_3,e_5\textgreater$
&
$\textbf{A}_{20}$
&
$\begin{array}{ll}
e_1 \cdot e_1=e_2\\
e_1 \cdot e_2=e_3\\
e_1 \cdot e_4=e_5
\end{array}$
$\begin{array}{ll}  
e_2 \cdot e_1=e_3\\
e_4 \cdot e_4=e_3+e_5
\end{array}$
&
$\textless e_2,e_3,e_5\textgreater$
\\ \hline 
$\textbf{A}_{21}$&
$\begin{array}{ll}  
e_1 \cdot e_1=e_2\\
e_1 \cdot e_2=e_3\\
e_1 \cdot e_4=e_5\\
\end{array}$
$\begin{array}{ll}  
e_2 \cdot e_1=e_3\\
e_4 \cdot e_1=e_2-e_5\\
e_5 \cdot e_1=e_3
\end{array}$
&
$\textless e_2,e_3\textgreater$
&
$\textbf{A}_{22}$
&
$\begin{array}{ll}
e_1 \cdot e_1=e_2\\
e_1 \cdot e_2=e_3\\
e_1 \cdot e_4=e_5\\
e_2 \cdot e_1=e_3\\
\end{array}$
$\begin{array}{ll}  
e_4 \cdot e_1=e_2-e_3\\
e_4 \cdot e_4=e_3\\
e_5 \cdot e_1=e_3
\end{array}$
&
$\textless e_2,e_3\textgreater$
\\ \hline 
$\textbf{A}_{23}$&
$\begin{array}{ll}  
e_1 \cdot e_1=e_2\\
e_1 \cdot e_2=e_3\\
e_1 \cdot e_4=e_5\\
e_2 \cdot e_1=e_3\\
\end{array}$
$\begin{array}{ll}  
e_4 \cdot e_1=e_2+e_5\\
e_4 \cdot e_2=2e_3\\
e_4 \cdot e_4=2e_3\\
e_5 \cdot e_1=e_3
\end{array}$
&
$\textless e_3\textgreater$
&
$\textbf{A}_{24}$
&
$\begin{array}{ll}
e_1 \cdot e_1=e_2\\
e_1 \cdot e_2=e_3\\
e_1 \cdot e_4=e_5\\
e_2 \cdot e_1=e_3\\
\end{array}$
$\begin{array}{ll}  
e_4 \cdot e_1=e_2+e_5\\
e_4 \cdot e_2=2e_3\\
e_4 \cdot e_4=e_3+2e_5\\
e_5 \cdot e_1=e_3
\end{array}$
&
$\textless e_3\textgreater$
\\ \hline 
$\textbf{A}_{25}$&
$\begin{array}{ll}  
e_1 \cdot e_1=e_2\\
e_1 \cdot e_2=e_3\\
e_1 \cdot e_4=e_5\\
e_2 \cdot e_1=e_3\\
\end{array}$
$\begin{array}{ll}  
e_4 \cdot e_1=e_3+e_5\\
e_4 \cdot e_4=2e_2\\
e_4 \cdot e_5=2e_3\\
e_5 \cdot e_4=e_3
\end{array}$
&
$\textless e_2, e_3\textgreater$
&
$\textbf{A}_{26}$
&
$\begin{array}{ll}
e_1 \cdot e_1=e_2\\
e_1 \cdot e_2=e_3\\
e_1 \cdot e_4=e_5\\
e_2 \cdot e_1=e_3\\
\end{array}$
$\begin{array}{ll}  
e_4 \cdot e_1=-e_5\\
e_4 \cdot e_4=2e_2\\
e_4 \cdot e_5=-e_3\\
e_5 \cdot e_4=e_3
\end{array}$
&
$\textless e_2, e_3\textgreater$
\\ \hline 
\end{tabular}
\end{center}

\small\begin{center}
    \begin{tabular}{||c||c||c||c||c||c||c||}
\hline
$\textbf{A}_i$&Multiplications&Center&$\textbf{A}_i$&Multiplications&Center.
\\ \hline 
$\textbf{A}_{27}$&
$\begin{array}{ll}  
e_1 \cdot e_1=e_2\\
e_1 \cdot e_2=e_3\\
e_1 \cdot e_4=e_5\\
e_2 \cdot e_1=e_3\\
\end{array}$
$\begin{array}{ll}  
e_4 \cdot e_1=e_2\\
e_4 \cdot e_2=2e_3\\
e_4 \cdot e_4=e_3+e_5\\
e_5 \cdot e_1=e_3
\end{array}$
&
$\textless e_3\textgreater$
&
$\textbf{A}_{28}$
&
$\begin{array}{ll}
e_1 \cdot e_1=e_2\\
e_1 \cdot e_2=e_3\\
e_1 \cdot e_4=e_5\\
e_2 \cdot e_1=e_3\\
\end{array}$
$\begin{array}{ll}  
e_4 \cdot e_1=e_3+e_5\\
e_4 \cdot e_4=-e_2+2e_5\\
e_4 \cdot e_5=e_3\\
e_5 \cdot e_4=-e_3
\end{array}$
&
$\textless e_2, e_3\textgreater$
\\ \hline 
$\textbf{A}_{29}$&
$\begin{array}{ll}  
e_1 \cdot e_1=e_2\\
e_1 \cdot e_2=e_3\\
e_1 \cdot e_4=e_5\\
e_2 \cdot e_1=e_3\\
\end{array}$
$\begin{array}{ll}
e_4 \cdot e_1=(1-i)e_2+ie_5\\
e_4 \cdot e_2=2e_3\\
e_4 \cdot e_4=-ie_2+e_3+(1+i)e_5\\
e_4 \cdot e_5=e_3\\
e_5 \cdot e_1=(1-i)e_3\\
e_5 \cdot e_4=-ie_3
\end{array}$
&
$\textless e_3\textgreater$
&
$\textbf{A}_{30}$
&
$\begin{array}{ll}
e_1 \cdot e_1=e_2\\
e_1 \cdot e_2=e_3\\
e_1 \cdot e_4=e_5\\
e_2 \cdot e_1=e_3\\
\end{array}$
$\begin{array}{ll} 
e_4 \cdot e_1=(1+i)e_2-ie_5\\
e_4 \cdot e_2=2e_3\\
e_4 \cdot e_4=ie_2+e_3+(1-i)e_5\\
e_4 \cdot e_5=e_3\\
e_5 \cdot e_1=(1+i)e_3\\
e_5 \cdot e_4=ie_3
\end{array}$
&
$\textless e_3\textgreater$
\\ \hline 
\end{tabular}
\end{center}
\begin{center}
\begin{tabular}{||c||c||c||c||c||c||c||}
\hline
$\textbf{A}_i$&Multiplications&Center.
\\ \hline  
$\textbf{A}_{31}$&
$\begin{array}{ll}  
e_1 \cdot e_1=e_2\\
e_1 \cdot e_2=e_3\\
e_1 \cdot e_4=e_5\\
e_2 \cdot e_1=e_3\\
e_4 \cdot e_1=(1-\alpha)e_2+\alpha e_5\\
e_4 \cdot e_2=(1-\alpha^2)e_3\\
e_4 \cdot e_4=-\alpha e_2+(1+\alpha)e_5\\
\end{array}$
$\begin{array}{ll}  
e_4 \cdot e_5=-\alpha^2e_3\\
e_5 \cdot e_1=(1-\alpha)e_3\\
e_5 \cdot e_4=-\alpha e_3
\end{array}$
&
$\textless e_3\textgreater$
\\ \hline   
\end{tabular}
\end{center}

\section{Derivation}
Let $(\textbf{A}, \cdot)$   be an $n$-dimensional Leibniz algebra with basis $\{e_i\}\, (1 \leq i \leq n)$ and let $d$ be a derivation on $\textbf{A}.$
For any $i,j,k\in\mathbb{N}$, $1 \leq i,j,k \leq n$, let us put
$$e_i\cdot e_j=\sum_{k=1}^n\lambda_{ij}^ke_k\quad \text{and}\quad d(e_i)=\sum_{j=1}^nd_{ji}e_j.$$
Then, in term of basis elements, equation (\ref{d2})  is equivalent to
\begin{align*}
\sum_{k=1}^n\Bigg(\lambda^k_{ij} d_{tk} -\lambda^t_{kj} d_{ki}-\lambda^t_{ik}d_{kj}\Bigg)=0,
\end{align*} for $i,j,t=1,2,\dots,n$.
\begin{proposition}
The description of the Derivation of every 5-dimensional  associative algebra is given below.
\end{proposition}
\begin{center}
\begin{tabular}{||c||c||c||c||c||c||c||}
\hline
$\textbf{A}_{i}$&\textbf{Der(A)}&\textbf{dim}&$\textbf{A}_{i}$&\textbf{Der(A)}&\textbf{dim}
\\ \hline 
$\textbf{A}_{1}$&
$\left(\begin{array}{cccccc}
0&0&0&0&0\\
d_{21}&d_{23}&0&0\\
-d_{23}&0&0&0&0\\
d_{41}&2d_{22}-d_{23}&d_{43}&0&d_{45}\\
d_{51}&0&d_{53}&0&d_{55}\\
\end{array}\right)$
&
8
&
$\textbf{A}_{2}$
&
$\left(\begin{array}{cccccc}
d_{11}&0&0&0&0\\
d_{21}&2d_{11}&0&0&0\\
d_{31}&0&2d_{11}&0&0\\
d_{41}&2d_{21}+d_{31}&d_{21}&3d_{11}&d_{45}\\
d_{51}&0&d_{53}&0&d_{55}\\
\end{array}\right)$
&
$8$
\\ \hline
$\textbf{A}_{3}$&
$\left(\begin{array}{cccccc}
d_{11}&0&0&0&0\\
d_{21}&2d_{11}&0&0&0\\
d_{31}&2d_{31}&3d_{11}&0&0\\
d_{41}&2d_{42}&2d_{21}&4d_{11}&d_{45}\\
d_{42}-2d_{31}&0&0&0&3d_{11}\\
\end{array}\right)$
&
6
&
$\textbf{A}_{4}$
&
$\left(\begin{array}{cccccc}
0&0&0&0&0\\
d_{21}&0&0&0&0\\
d_{31}&2d_{31}&0&0&d_{35}\\
d_{41}&2d_{31}-d_{35}&3d_{21}&0&d_{45}\\
-d_{35}&0&0&0&0\\
\end{array}\right)$
&
$4$
\\ \hline
$\textbf{A}_{5}$&
$\left(\begin{array}{cccccc}
\frac{d_{33}}{3}&0&0&0&0\\
d_{21}&\frac{2d_{33}}{3}&0&d_{24}&d_{25}\\
d_{31}&2d_{21}-d_{25}&d_{33}&d_{34}&d_{35}\\
-d_{25}&0&0&\frac{2d_{33}}{3}&0\\
-d_{24}&0&0&0&\frac{d_{33}}{2}
\end{array}\right)$
&
8
&
$\textbf{A}_{6}$
&
$\left(\begin{array}{cccccc}
\frac{d_{33}}{3}&0&0&0&0\\
d_{21}&\frac{2d_{33}}{3}&0&0&d_{25}\\
d_{31}&d_{32}&d_{33}&d_{34}&d_{35}\\
d_{32}-2d_{21}&0&0&\frac{2d_{33}}{3}&0\\
-d_{25}&0&0&0&\frac{2d_{33}}{3}
\end{array}\right)$
&
$8$
\\ \hline 
$\textbf{A}_{7}$&
$\left(\begin{array}{cccccc}
0&0&0&0&0\\
d_{21}&0&0&d_{24}&d_{25}\\
d_{31}&2d_{21}-d_{24}&0&d_{34}&d_{35}\\
-d_{24}&0&0&0&0\\
-d_{25}&0&0&0&0
\end{array}\right)$
&
6
&
$\textbf{A}_{8}$
&
$\left(\begin{array}{cccccc}
\frac{d_{33}}{3}&0&0&0&0\\
d_{21}&\frac{2d_{33}}{3}&0&0&0\\
d_{31}&2d_{21}&d_{33}&d_{34}&d_{35}\\
0&0&0&\frac{d_{33}}{3}&d_{45}\\
0&0&0&0&\frac{d_{33}}{3}
\end{array}\right)$
&
$6$
\\ \hline 
\end{tabular}
\end{center}	

\begin{center}
\begin{tabular}{||c||c||c||c||c||c||c||}
\hline
$\textbf{A}_{i}$&\textbf{Der(A)}&\textbf{dim}&$\textbf{A}_{i}$&\textbf{Der(A)}&\textbf{dim}
\\ \hline  
$\textbf{A}_{9}$&
$\left(\begin{array}{cccccc}
\frac{d_{33}}{3}&0&0&0&0\\
d_{21}&\frac{2d_{33}}{3}&0&0&0\\
d_{31}&2d_{21}&d_{33}&d_{34}&d_{35}\\
0&0&0&0&0\\
0&0&0&d_{44}&d_{33}-d_{44}
\end{array}\right)$
&
7
&
$\textbf{A}_{10}$
&
$\left(\begin{array}{cccccc}
d_{11}&0&0&0&0\\
d_{21}&2d_{11}&0&0&0\\
d_{31}&2d_{21}&3d_{11}&d_{34}&0\\
d_{41}&0&0&d_{44}&0\\
d_{51}&d_{41}&0&d_{54}&d_{11}+d_{44}
\end{array}\right)$
&
$8$
\\ \hline 
$\textbf{A}_{11}$&
$\left(\begin{array}{cccccc}
\frac{d_{33}}{3}&0&0&0&0\\
d_{21}&\frac{2d_{33}}{3}&0&0&0\\
d_{31}&2d_{21}&d_{33}&d_{34}&0\\
-d_{24}&0&0&\frac{d_{33}}{2}&0\\
d_{51}&-d_{21}&0&d_{54}&\frac{5d_{33}}{6}
\end{array}\right)$
&
7
&
$\textbf{A}_{12}$
&
$\small\left(\begin{array}{cccccc}
0&0&0&0&0\\
d_{21}&d_{22}&0&0&0\\
d_{31}&d_{23}&d_{33}&d_{34}&d_{35}\\
k&0&0&d_{33}-d_{22}&0\\
d_{32}-d_{21}&-k&0&0&d_{33}-\frac{d_{22}}{2}
\end{array}\right)$
&
$6$
\\ \hline 
$\textbf{A}_{13}$&
$\left(\begin{array}{cccccc}
0&0&0&0&0\\
d_{21}&0&0&0&0\\
d_{31}&d_{32}&0&d_{34}&d_{35}\\
0&0&0&0&0\\
d_{32}-2d_{21}&0&0&-d_{21}&0
\end{array}\right)$
&
7
&
$\textbf{A}_{14}$
&
$\small\left(\begin{array}{cccccc}
0&0&0&0&0\\
d_{21}&d_{22}&0&0&0\\
d_{31}&d_{23}&d_{33}&d_{34}&d_{35}\\
k&0&0&d_{33}-d_{22}&0\\
d_{32}-d_{21}&-k&0&0&d_{33}-\frac{d_{22}}{2}
\end{array}\right)$
&
$6$
\\ \hline 
$\textbf{A}_{15}$&
$\left(\begin{array}{cccccc}
d_{11}&0&0&0&0\\
d_{21}&2d_{11}&0&2d_{24}&0\\
d_{31}&2d_{21}-d_{24}&3d_{11}&d_{34}&\frac{d_{11}}{2}\\
-d_{24}&0&0&\frac{3d_{11}}{2}&0\\
d_{51}&-d_{24}&0&d_{54}&\frac{5d_{11}}{2}
\end{array}\right)$
&
7
&
$\textbf{A}_{16}$
&
$\left(\begin{array}{cccccc}
d_{11}&0&0&0&0\\
d_{21}&d_{11}&0&0&0\\
d_{31}&d_{23}&d_{11}&d_{34}&0\\
0&0&0&d_{11}&0\\
d_{51}&0&0&d_{54}&d_{11}
\end{array}\right)$
&
$6$
\\ \hline 
$\textbf{A}_{17}$&
$\left(\begin{array}{cccccc}
d_{11}&0&0&0&0\\
d_{21}&2d_{11}&0&-d_{41}&0\\
d_{31}&2d_{21}&3d_{11}&d_{34}&0\\
d_{41}&0&0&\frac{3d_{11}}{2}&0\\
d_{51}&(\alpha-1)d_{41}&0&d_{54}&\frac{5d_{11}}{2}
\end{array}\right)$
&
7
&
$\textbf{A}_{18}$
&
$\left(\begin{array}{cccccc}
0&0&0&0&0\\
d_{21}&0&0&0&0\\
d_{31}&2d_{23}&0&d_{34}&0\\
0&0&0&0&0\\
d_{51}&0&0&d_{54}&0
\end{array}\right)$
&
$6$
\\ \hline 
$\textbf{A}_{19}$&
$\left(\begin{array}{cccccc}
d_{11}&0&0&0&0\\
d_{21}&2d_{11}&0&0\\
d_{31}&2d_{21}&2d_{11}&d_{34}&0\\
0&0&0&2d_{11}&0\\
d_{51}&0&0&d_{54}&4d_{11}\\
\end{array}\right)$
&
6
&
$\textbf{A}_{20}$
&
$\left(\begin{array}{cccccc}
d_{11}&0&0&0&0\\
d_{21}&2d_{11}&0&0&0\\
d_{31}&2d_{21}&3d_{11}&d_{34}&d_{11}\\
0&0&0&2d_{11}&0\\
d_{51}&0&0&d_{54}&4d_{11}\\
\end{array}\right)$
&
$6$
\\ \hline 
$\textbf{A}_{21}$&
$\left(\begin{array}{cccccc}
d_{11}&0&0&0&0\\
d_{21}&d_{22}&0&d_{24}&0\\
d_{31}&d_{32}&d_{11}+d_{22}&d_{34}&d_{24}\\
d_{22}-2d_{24}&0&0&\frac{3d_{11}}{2}&0\\
k_1&0&0&\frac{-d_{11}}{2}&2d_{11}\\
\end{array}\right)$
&
6
&
$\textbf{A}_{22}$
&
$\left(\begin{array}{cccccc}
d_{11}&0&0&0&0\\
d_{21}&d_{11}&0&d_{24}&0\\
d_{31}&d_{32}&3d_{11}&d_{34}&d_{35}\\
0&2d_{31}-d_{35}&3d_{21}&0&d_{45}\\
-d_{35}&0&0&0&0\\
\end{array}\right)$
&
$4$
\\ \hline 
$\textbf{A}_{23}$&
$\left(\begin{array}{cccccc}
d_{11}&0&0&0&0\\
d_{21}&2d_{11}&0&d_{24}&0\\
d_{31}&d_{32}&2d_{11}&d_{34}&d_{24}\\
0&0&0&d_{11}&0\\
k_2&0&0&k_2&2d_{11}
\end{array}\right)$
&
7
&
$\textbf{A}_{24}$
&
$\left(\begin{array}{cccccc}
0&0&0&0&0\\
d_{21}&0&0&d_{24}&0\\
d_{31}&d_{32}&0&d_{34}&d_{24}\\
0&0&0&0&0\\
d_{32}-2d_{21}&0&0&d_{32}-2d_{21}&0
\end{array}\right)$
&
$5$
\\ \hline 
$\textbf{A}_{25}$&
$\left(\begin{array}{cccccc}
0&0&0&0&0\\
d_{21}&0&0&d_{24}&0\\
d_{31}&2d_{21}&0&d_{34}&d_{35}\\
0&0&0&0&0\\
d_{35}-d_{25}&0&0&d_{21}&0
\end{array}\right)$
&
6
&
$\textbf{A}_{26}$
&
$\left(\begin{array}{cccccc}
d_{11}&0&0&0&0\\
0&2d_{11}&0&0&0\\
d_{31}&0&3d_{11}&d_{34}&d_{35}\\
0&0&0&d_{11}&0\\
d_{35}&0&0&d_{54}&d_{11}
\end{array}\right)$
&
$5$
\\ \hline 
$\textbf{A}_{27}$&
$\left(\begin{array}{cccccc}
d_{11}&0&0&0&0\\
d_{21}&d_{11}&0&d_{24}&0\\
d_{31}&d_{32}&d_{11}&d_{34}&d_{24}-d_{11}\\
-d_{11}&0&0&0&0\\
k_4&0&0&k_5&0
\end{array}\right)$
&
6
&
$\textbf{A}_{28}$
&
$\left(\begin{array}{cccccc}
d_{11}&0&0&0&0\\
d_{21}&2d_{11}&0&\frac{d_{11}}{2}+d_{21}&0\\
d_{31}&2d_{21}&3d_{11}&d_{34}&d_{21}\\
0&0&0&d_{11}&0\\
\frac{d_{11}}{2}&0&0&d_{54}&2d_{11}
\end{array}\right)$
&
$8$
\\ \hline 
$\textbf{A}_{29}$&
$\left(\begin{array}{cccccc}
0&0&0&0&0\\
0&0&0&0&0\\
d_{31}&0&0&d_{34}&0\\
0&0&0&0&0\\
0&0&0&0&0\\
\end{array}\right)$
&
2
&
$\textbf{A}_{30}$
&
$\left(\begin{array}{cccccc}
0&0&0&0&0\\
0&0&0&0&0\\
d_{31}&0&0&d_{34}&0\\
0&0&0&0&0\\
0&0&0&0&0\\
\end{array}\right)$
&
$6$
\\ \hline 
\end{tabular}
\end{center}	
\small
\begin{center}
\begin{tabular}{||c||c||c||c||c||c||c||}
\hline
$\textbf{A}_{i}$&\textbf{Der(D)}&\textbf{dim}
\\ \hline
$\textbf{A}_{31}$&
$\left(\begin{array}{cccccc}
(\alpha^2+\alpha-1)d_{41}+d_{44}&0&0&(\alpha^2+\alpha)d_{41}&0\\
-d_{51}&(\alpha^2+\alpha-1)d_{41}+2d_{44}&0&-d_{54}&\alpha^2d_{41}\\
d_{31}&-(\alpha+1)d_{51}&(2\alpha^2+2\alpha-1)d_{41}+3d_{44}&-\alpha d_{51}+d_{54}&d_{35}\\
d_{41}&0&0&d_{44}&0\\
d_{51}&(\alpha+1)d_{41}&0&d_{54}&\alpha(2+\alpha)d_{41}+2d_{44}
\end{array}\right)$
&
$
6
$
\\ \hline 
\end{tabular}
\end{center}	
\section{Element of centroids}
Let $(\textbf{A},\cdot)$   be an $n$-dimensional associative algebra with basis $\{e_i\}\, (1 \leq i \leq n)$ and let $c$ be the Centroid on $\textbf{A}$.
For any $i,j,k\in\mathbb{N}$, $1 \leq i,j,k \leq n$, let us put
$$e_i \cdot  e_j=\sum_{k=1}^n\lambda_{ij}^ke_k\quad \text{and}\quad c(e_i)=\sum_{j=1}^nc_{ji}e_j.$$
Then, in term of basis elements, equation (\ref{d2}) is equivalent to
\begin{align*} 
   \sum^{n}_{k=1}c_{ki}\lambda_{kj}^{p}=\sum^{n}_{k=1}\lambda_{ij}^{k}c_{pk}=\sum^{n}_{k=1}c_{kj}\lambda_{ik}^{p},
\end{align*} for $i,j,p=1,2,\dots,n$.

\begin{proposition}
The description of the Derivation of every 5-dimensional  associative algebra is given below.
\end{proposition}
\begin{center}
\begin{tabular}{||c||c||c||c||c||c||c||}
\hline
$\textbf{A}_{i}$&\textbf{Cent(A)}&\textbf{dim}&$\textbf{A}_{i}$&\textbf{Cent(A)}&\textbf{dim}
\\ \hline 
$\textbf{A}_{1}$&
$\left(\begin{array}{cccccc}
c_{33}&0&0&0&0\\
c_{21}&c_{33}&0&0&0\\
0&0&c_{33}&0&0\\
c_{41}&c_{21}&c_{43}&c_{33}&c_{45}\\
c_{51}&0&c_{53}&0&c_{55}
\end{array}\right)$
&
8
&
$\textbf{A}_{2}$
&
$\left(\begin{array}{cccccc}
c_{11}&0&0&0&0\\
c_{21}&c_{11}&0&0&0\\
0&0&c_{11}&0&0\\
c_{41}&c_{21}&c_{43}&c_{11}&c_{45}\\
c_{51}&0&c_{53}&0&c_{55}
\end{array}\right)$
&
$8$
\\ \hline 
$\textbf{A}_{3}$&
$\left(\begin{array}{cccccc}
c_{11}&0&0&0&0\\
c_{21}&c_{11}&0&0&0\\
0&0&c_{11}&0&0\\
c_{41}&c_{21}&c_{43}&c_{11}&c_{45}\\
c_{51}&0&c_{53}&0&c_{11}
\end{array}\right)$
&
7
&
$\textbf{A}_{4}$
&
$\left(\begin{array}{cccccc}
c_{44}&0&0&0&0\\
c_{32}&c_{44}&0&0&0\\
c_{31}&c_{32}&c_{44}&0&0\\
c_{41}&c_{31}&c_{32}&c_{44}&c_{45}\\
0&0&0&0&c_{44}
\end{array}\right)$
&
$5$
\\ \hline 
$\textbf{A}_{5}$&
$\left(\begin{array}{cccccc}
c_{33}&0&0&0&0\\
c_{21}&c_{33}&0&0&0\\
c_{31}&c_{21}&c_{33}&c_{34}&0\\
0&0&0&c_{33}&c_{35}\\
0&0&0&0&c_{33}
\end{array}\right)$
&
5
&
$\textbf{A}_{6}$
&
$\left(\begin{array}{cccccc}
c_{33}&0&0&0&0\\
c_{21}&c_{33}&0&0&0\\
c_{31}&c_{21}&c_{33}&c_{34}&c_{35}\\
0&0&0&c_{33}&0\\
0&0&0&0&c_{33}
\end{array}\right)$
&
$5$
\\ \hline 
$\textbf{A}_{7}$&
$\left(\begin{array}{cccccc}
c_{33}&0&0&0&0\\
c_{21}&c_{33}&0&0&0\\
c_{31}&c_{21}&c_{33}&c_{34}&c_{35}\\
0&0&0&c_{33}&0\\
0&0&0&0&c_{33}
\end{array}\right)$
&
5
&
$\textbf{A}_{8}$
&
$\left(\begin{array}{cccccc}
c_{33}&0&0&0&0\\
c_{21}&c_{33}&0&0&0\\
c_{31}&c_{21}&c_{33}&c_{34}&c_{35}\\
0&0&0&c_{33}&0\\
0&0&0&0&c_{33}
\end{array}\right)$
&
$5$
\\ \hline 
$\textbf{A}_{9}$&
$\left(\begin{array}{cccccc}
c_{33}&0&0&0&0\\
c_{21}&c_{33}&0&0&0\\
c_{31}&c_{21}&c_{33}&c_{34}&c_{35}\\
0&0&0&c_{33}&0\\
0&0&0&0&c_{33}
\end{array}\right)$
&
5
&
$\textbf{A}_{10}$
&
$\left(\begin{array}{cccccc}
c_{11}&0&0&0&0\\
c_{21}&c_{11}&0&0&0\\
c_{31}&c_{21}&c_{11}&c_{34}&0\\
0&0&0&c_{11}&0\\
c_{51}&0&0&c_{54}&c_{11}
\end{array}\right)$
&
$6$
\\ \hline 
$\textbf{A}_{11}$&
$\left(\begin{array}{cccccc}
c_{33}&0&0&0&0\\
c_{21}&c_{33}&0&0&0\\
c_{31}&c_{21}&c_{33}&c_{34}&0\\
0&0&0&c_{33}&0\\
c_{51}&0&0&c_{54}&c_{33}
\end{array}\right)$
&
6
&
$\textbf{A}_{12}$
&
$\left(\begin{array}{cccccc}
c_{33}&0&0&0&0\\
0&c_{33}&0&0&0\\
c_{31}&0&c_{33}&c_{34}&c_{35}\\
0&0&0&c_{33}&0\\
0&0&0&0&c_{33}
\end{array}\right)$
&
$4$
\\ \hline 
\end{tabular}
\end{center}	

\begin{center}
\begin{tabular}{||c||c||c||c||c||c||c||}
\hline
$\textbf{A}_{i}$&\textbf{Cent(A)}&\textbf{dim}&$\textbf{A}_{i}$&\textbf{Cent(A)}&\textbf{dim}
\\ \hline 
$\textbf{A}_{13}$&
$\left(\begin{array}{cccccc}
c_{33}&0&0&0&0\\
0&0&0&0\\
c_{31}&0&c_{33}&c_{34}&c_{35}\\
0&0&0&c_{33}&0\\
0&0&0&0&c_{33}
\end{array}\right)$
&
4
&
$\textbf{A}_{14}$
&
$\left(\begin{array}{cccccc}
c_{11}&0&0&0&0\\
c_{21}&c_{11}&0&0&0\\
c_{31}&c_{21}&c_{11}&c_{34}&0\\
0&0&0&c_{33}&0\\
c_{51}&0&0&c_{54}&c_{33}
\end{array}\right)$
&
$6$
\\ \hline 
$\textbf{A}_{15}$&
$\left(\begin{array}{cccccc}
c_{11}&0&0&0&0\\
c_{21}&c_{11}&0&0&0\\
c_{31}&c_{21}&c_{11}&c_{34}&0\\
0&0&0&c_{11}&0\\
c_{51}&0&0&c_{54}&c_{11}
\end{array}\right)$
&
6
&
$\textbf{A}_{16}$
&
$\left(\begin{array}{cccccc}
c_{55}&0&0&0&0\\
c_{21}&c_{55}&0&0&0\\
c_{31}&c_{21}&c_{55}&c_{34}&0\\
0&0&0&c_{55}&0\\
c_{51}&0&0&c_{54}&c_{55}
\end{array}\right)$
&
$6$
\\ \hline 
$\textbf{A}_{17}$&
$\left(\begin{array}{cccccc}
c_{55}&0&0&0&0\\
c_{21}&c_{55}&0&0&0\\
c_{31}&c_{21}&c_{55}&c_{34}&0\\
0&0&0&c_{55}&0\\
c_{51}&0&0&c_{54}&c_{55}
\end{array}\right)$
&
6
&
$\textbf{A}_{18}$
&
$\left(\begin{array}{cccccc}
c_{11}&0&0&0&0\\
c_{21}&c_{11}&0&0&0\\
c_{31}&c_{21}&c_{11}&c_{34}&0\\
0&0&0&c_{11}&0\\
c_{51}&0&0&c_{54}&c_{11}
\end{array}\right)$
&
$6$
\\ \hline
$\textbf{A}_{19}$&
$\left(\begin{array}{cccccc}
c_{11}&0&0&0&0\\
c_{21}&c_{11}&0&0&0\\
c_{31}&c_{21}&c_{11}&c_{34}&0\\
0&0&0&c_{11}&0\\
c_{51}&0&0&c_{54}&c_{11}
\end{array}\right)$
&
6
&
$\textbf{A}_{20}$
&
$\left(\begin{array}{cccccc}
c_{11}&0&0&0&0\\
c_{21}&c_{11}&0&0&0\\
c_{31}&c_{21}&c_{11}&c_{34}&0\\
0&0&0&c_{11}&0\\
c_{51}&0&0&c_{54}&c_{11}
\end{array}\right)$
&
$6$
\\ \hline 
$\textbf{A}_{21}$&
$\left(\begin{array}{cccccc}
c_{11}&0&0&0&0\\
0&c_{11}&0&0&0\\
c_{31}&c_{21}&c_{11}&c_{34}&0\\
0&0&0&c_{11}&0\\
0&0&0&0&c_{11}
\end{array}\right)$
&
4
&
$\textbf{A}_{22}$
&
$\left(\begin{array}{cccccc}
c_{11}&0&0&0&0\\
c_{21}&c_{11}&0&0&0\\
c_{31}&c_{21}&c_{11}&c_{34}&c_{35}\\
0&0&0&c_{11}&0\\
0&0&0&c_{54}&c_{11}
\end{array}\right)$
&
$6$
\\ \hline 
$\textbf{A}_{23}$&
$\left(\begin{array}{cccccc}
c_{11}&0&0&0&0\\
0&c_{11}&0&0&0\\
c_{31}&c_{21}&c_{11}&c_{34}&0\\
0&0&0&c_{11}&0\\
0&0&0&0&c_{11}
\end{array}\right)$
&
4
&
$\textbf{A}_{24}$
&
$\left(\begin{array}{cccccc}
c_{11}&0&0&0&0\\
0&c_{11}&0&0&0\\
c_{31}&c_{21}&c_{11}&c_{34}&0\\
0&0&0&c_{11}&0\\
0&0&0&0&c_{11}
\end{array}\right)$
&
$4$
\\ \hline
$\textbf{A}_{25}$&
$\left(\begin{array}{cccccc}
c_{11}&0&0&0&0\\
c_{21}&c_{11}&0&c_{24}&0\\
c_{31}&c_{21}&c_{11}&c_{34}&c_{24}\\
0&0&0&c_{11}&0\\
c_{24}&0&0&c_{21}&c_{11}
\end{array}\right)$
&
7
&
$\textbf{A}_{26}$
&
$\left(\begin{array}{cccccc}
c_{11}&0&0&0&0\\
0&c_{11}&0&0&0\\
c_{31}&c_{21}&c_{11}&c_{34}&0\\
0&0&0&c_{11}&0\\
0&0&0&0&c_{11}
\end{array}\right)$
&
$3$
\\ \hline 
$\textbf{A}_{27}$&
$\left(\begin{array}{cccccc}
c_{11}&0&0&0&0\\
c_{21}&c_{11}&0&0&0\\
c_{31}&c_{21}&c_{11}&c_{34}&0\\
0&0&0&c_{11}&0\\
0&0&0&c_{21}&c_{11}
\end{array}\right)$
&
4
&
$\textbf{A}_{28}$
&
$\left(\begin{array}{cccccc}
c_{11}&0&0&0&0\\
0&c_{11}&0&0&0\\
c_{31}&c_{21}&c_{11}&c_{34}&0\\
0&0&0&c_{11}&0\\
0&0&0&0&c_{11}
\end{array}\right)$
&
$4$
\\ \hline 
$\textbf{A}_{29}$&
$\left(\begin{array}{cccccc}
c_{44}&0&0&0&0\\
0&c_{44}&0&0&0\\
c_{31}&0&c_{44}&0&0\\
0&0&0&c_{44}&0\\
0&0&0&0&c_{44}
\end{array}\right)$
&
2
&
$\textbf{A}_{30}$
&
$\left(\begin{array}{cccccc}
c_{44}&0&0&0&0\\
0&c_{44}&0&0&0\\
c_{31}&0&c_{44}&0&0\\
0&0&0&c_{44}&0\\
0&0&0&0&c_{44}
\end{array}\right)$
&
$2$
\\ \hline 
$\textbf{A}_{31}$&
$\left(\begin{array}{cccccc}
c_{44}&0&0&0&0\\
0&c_{44}&0&0&0\\
c_{31}&0&c_{44}&c_{34}&0\\
0&0&0&c_{44}&0\\
0&0&0&0&c_{44}
\end{array}\right)$
&
3
&
&
&
\\ \hline 
\end{tabular}
\end{center}

\section{Inner Derivation}
Let $(\textbf{A},\cdot)$  be an  $n$-dimensional associative algebra. Fix a basis  $\{e_1, e_2, \dots, e_n\}$ of  $\textbf{A}$  and a vector  
$w = a_1 e_1 + a_2 e_2 + \dots + a_n e_n$ in  $\textbf{A}$. Let $ad_w$  be an inner derivation of  $\textbf{A}$. In particular, we have:

\begin{align*}
ad_w(e_i) &= e_i\circ w - w\circ e_i, \quad \forall i \in \{1, 2, \dots, n\}.
\end{align*}
Moreover, $ad_w$ can be considered as a linear transformation of  $A$, represented in a matrix form:
\begin{align*}
ad_w(e_i) &= \sum_{j=1}^n d_{ij} e_j, \quad i = 1, 2, \dots, n.
\end{align*}
This leads to:
\begin{align*}
\sum_{j=1}^n d_{ij} e_j &= e_i\circ w - w\circ e_i, \quad i = 1, 2, \dots, n.
\end{align*}
Then, we obtain the following system of equations for the coefficients  $d_{ij}$:
\begin{align*}
d_{ij} &= \sum_{t=1}^n\Bigg(a_tC_{it}^j-a_tC_{ti}^j\Bigg), \quad i, j = 1, 2, \dots, n. 
\end{align*}
Solving this system of equation, we find the inner derivation in matrix form.
\begin{proposition}
The description of the Derivation of every 5-dimensional  associative algebra is given below.
\end{proposition}
\begin{center}
\begin{tabular}{||c||c||c||c||c||c||c||}
\hline
$\textbf{A}_{i}$&$\textbf{Inn(A)}$&\textbf{dim}&$\textbf{A}_{i}$&$\textbf{Inn(A)}$&\textbf{dim}
\\ \hline 
$\textbf{A}_{1}$&
$\left(\begin{array}{cccccc}
0&0&0&0&0\\
0&0&0&0&0\\
0&0&0&0&0\\
a_{3}&0&-a_{1}&0&0\\
0&0&0&0&0
\end{array}\right)$
&
2
&
$\textbf{A}_{2}$
&
$\left(\begin{array}{cccccc}
0&0&0&0&0\\
0&0&0&0&0\\
0&0&0&0&0\\
a_{3}&0&-a_{1}&0&0\\
0&0&0&0&0
\end{array}\right)$
&
$2$
\\ \hline 
$\textbf{A}_{3}$&
$\left(\begin{array}{cccccc}
0&0&0&0&0\\
0&0&0&0&0\\
0&0&0&0&0\\
a_{5}&0&0&0&-a_1\\
0&0&0&0&0
\end{array}\right)$
&
6
&
$\textbf{A}_{4}$
&
$\left(\begin{array}{cccccc}
0&0&0&0&0\\
0&0&0&0&0\\
0&0&0&0&0\\
a_{5}&0&0&0&-a_1\\
0&0&0&0&0
\end{array}\right)$
&
$2$
\\ \hline 
$\textbf{A}_{5}$&
$\left(\begin{array}{cccccc}
0&0&0&0&0\\
0&0&0&0&0\\
a_4&0&0&-a_1&0\\
0&0&0&0&0\\
0&0&0&0&0
\end{array}\right)$
&
2
&
$\textbf{A}_{6}$
&
$\left(\begin{array}{cccccc}
0&0&0&0&0\\
0&0&0&0&0\\
a_4&0&0&-a_1&0\\
0&0&0&0&0\\
0&0&0&0&0
\end{array}\right)$
&
$2$
\\ \hline 
$\textbf{A}_{7}$&
$\left(\begin{array}{cccccc}
0&0&0&0&0\\
0&0&0&0&0\\
a_4&0&0&-a_1&0\\
0&0&0&0&0\\
0&0&0&0&0
\end{array}\right)$
&
2
&
$\textbf{A}_{8}$
&
$\left(\begin{array}{cccccc}
0&0&0&0&0\\
0&0&0&0&0\\
0&0&0&2a_5&-2a_4\\
0&0&0&0&0\\
0&0&0&0&0
\end{array}\right)$
&
$2$
\\ \hline 
$\textbf{A}_{9}$&
$\left(\begin{array}{cccccc}
0&0&0&0&0\\
0&0&0&0&0\\
0&0&0&(1-\alpha)a_5&(\alpha-1)a_4\\
0&0&0&0&0\\
0&0&0&0&0
\end{array}\right)$
&
2
&
$\textbf{A}_{10}$
&
$\left(\begin{array}{cccccc}
0&0&0&0&0\\
0&0&0&0&0\\
0&0&0&0&0\\
0&0&0&0&0\\
-a_4&0&0&a_1&0
\end{array}\right)$
&
$2$
\\ \hline 
$\textbf{A}_{11}$&
$\left(\begin{array}{cccccc}
0&0&0&0&0\\
0&0&0&0&0\\
0&0&0&0&0\\
0&0&0&0&0\\
-a_4&0&0&a_1&0
\end{array}\right)$
&
6
&
$\textbf{A}_{12}$
&
$\left(\begin{array}{cccccc}
0&0&0&0&0\\
0&0&0&0&0\\
-a_5&-a_4&0&a_2&a_1\\
0&0&0&0&0\\
-a_4&0&0&a_1&0
\end{array}\right)$
&
$4$
\\ \hline 
$\textbf{A}_{13}$&
$\left(\begin{array}{cccccc}
0&0&0&0&0\\
0&0&0&0&0\\
-a_4-a_5&-a_4&0&a_1+a_2&a_1\\
0&0&0&0&0\\
0&0&0&0&0
\end{array}\right)$
&
4
&
$\textbf{A}_{14}$
&
$\left(\begin{array}{cccccc}
0&0&0&0&0\\
0&0&0&0&0\\
-a_4&0&0&a_1&0\\
0&0&0&0&0\\
0&0&0&0&0
\end{array}\right)$
&
$2$
\\ \hline 
$\textbf{A}_{15}$&
$\left(\begin{array}{cccccc}
0&0&0&0&0\\
0&0&0&0&0\\
-a_4&0&0&a_1&0\\
0&0&0&0&0\\
-a_5&0&0&a_1&0
\end{array}\right)$
&
3
&
$\textbf{A}_{16}$
&
$\left(\begin{array}{cccccc}
0&0&0&0&0\\
0&0&0&0&0\\
0&0&0&0&0\\
0&0&0&0&0\\
(1-\alpha)a_4&0&0&(\alpha-1)a_1&0
\end{array}\right)$
&
$2$
\\ \hline 
\end{tabular}
\end{center}

\begin{center}
\begin{tabular}{||c||c||c||c||c||c||c||}
\hline
$\textbf{A}_{i}$&$\textbf{Inn(A)}$&\textbf{dim}&$\textbf{A}_{i}$&$\textbf{Inn(A)}$&\textbf{dim}
\\ \hline 
$\textbf{A}_{17}$&
$\left(\begin{array}{cccccc}
0&0&0&0&0\\
0&0&0&0&0\\
0&0&0&0&0\\
0&0&0&0&0\\
(1-\alpha)a_4&0&0&(\alpha-1)a_1&0
\end{array}\right)$
&
2
&
$\textbf{A}_{18}$
&
$\left(\begin{array}{cccccc}
0&0&0&0&0\\
0&0&0&0&0\\
-a_4&0&0&a_1&0\\
0&0&0&0&0\\
a_4&0&0&a_1&0
\end{array}\right)$
&
$2$
\\ \hline
$\textbf{A}_{19}$&
$\left(\begin{array}{cccccc}
0&0&0&0&0\\
0&0&0&0&0\\
-a_4&0&0&a_1&0\\
0&0&0&0&0\\
0&0&0&0&0
\end{array}\right)$
&
2
&
$\textbf{A}_{20}$
&
$\left(\begin{array}{cccccc}
0&0&0&0&0\\
0&0&0&0&0\\
-a_4&0&0&a_1&0\\
0&0&0&0&0\\
0&0&0&0&0
\end{array}\right)$
&
$2$
\\ \hline 
$\textbf{A}_{21}$&
$\left(\begin{array}{cccccc}
0&0&0&0&0\\
-a_4&0&0&a_1&0\\
-a_5&0&0&a_1&0\\
0&0&0&0&0\\
2a_4&0&0&-2a_1&0\\
\end{array}\right)$
&
3
&
$\textbf{A}_{22}$
&
$\left(\begin{array}{cccccc}
0&0&0&0&0\\
0&0&0&0&0\\
a_4-a_5&0&-a_1&a_1&0\\
0&0&0&0&0\\
0&0&0&0&0
\end{array}\right)$
&
$3$
\\ \hline 
$\textbf{A}_{23}$&
$\left(\begin{array}{cccccc}
0&0&0&0&0\\
0&0&0&a_1&0\\
-a_4&-2a_4&0&2a_2&a_1\\
0&0&0&0&0\\
-a_5&0&0&0&0
\end{array}\right)$
&
4
&
$\textbf{A}_{24}$
&
$\left(\begin{array}{cccccc}
0&0&0&0&0\\
-a_4&0&0&a_1&0\\
-a_5&-2a_4&0&2a_2&a_1\\
0&0&0&0&0\\
0&0&0&0&0
\end{array}\right)$
&
$4$
\\ \hline 
$\textbf{A}_{25}$&
$\left(\begin{array}{cccccc}
0&0&0&0&0\\
0&0&0&0&0\\
-a_4&0&0&a_1&0\\
0&0&0&0&0\\
0&0&0&0&0
\end{array}\right)$
&
2
&
$\textbf{A}_{26}$
&
$\left(\begin{array}{cccccc}
0&0&0&0&0\\
0&0&0&0&0\\
0&0&0&-2a_5&2a_4\\
0&0&0&0&0\\
2a_4&0&0&-2a_1&0
\end{array}\right)$
&
$3$
\\ \hline 
$\textbf{A}_{27}$&
$\left(\begin{array}{cccccc}
0&0&0&0&0\\
-a_4&0&0&0&a_1\\
-a_5&-a_4&0&0&a_2\\
0&0&0&0&0\\
a_4&0&0&-a_1&0
\end{array}\right)$
&
4
&
$\textbf{A}_{28}$
&
$\left(\begin{array}{cccccc}
0&0&0&0&0\\
0&0&0&0&0\\
-a_4&0&0&a_1-2a_5&-2a_4\\
0&0&0&0&0\\
0&0&0&0
\end{array}\right)$
&
$3$
\\ \hline 
$\textbf{A}_{29}$&
$\left(\begin{array}{cccccc}
0&0&0&0&0\\
0&-2a_4&0&2a_2-3a_1&0\\
0&0&0&5a_5&-4a_4-6a_4\\
0&0&0&0&0\\
0&0&0&3a_1&0
\end{array}\right)$
&
4
&
$\textbf{A}_{30}$
&
$\left(\begin{array}{cccccc}
0&0&0&0&0\\
-2a_4&-2a_4&0&5a_1+2a_2&0\\
-2a_5&0&0&-3a_5&6a_1+4a_4\\
0&0&0&0&0\\
2a_4&0&0&-5a_1&0
\end{array}\right)$
&
$4$
\\ \hline 
\end{tabular}
\end{center}	
\begin{center}
\begin{tabular}{||c||c||c||c||c||c||c||}
\hline
$\textbf{A}_{i}$&Inn(\textbf{A})&\textbf{dim}
\\ \hline
$\textbf{A}_{31}$&
$\left(\begin{array}{cccccc}
0&0&0&0&0\\
(\alpha-1)a_4&0&0&(\alpha-1)a_1&0\\
(\alpha-1)a_5&(\alpha^2-1)a_4&0&(1-\alpha^2)a_2-\alpha(1-\alpha)a_5&(1-\alpha)a_1-\alpha(1+\alpha)a_4\\
0&0&0&0&0\\
a_4-\alpha a_5&0&0&(\alpha-1)a_1&0
\end{array}\right)$
&
$
4
$
\\ \hline 
\end{tabular}
\end{center}	

\begin{corollary}\,
\begin{itemize}
    \item The dimensions of the derivation of associative algebras of $5$-dimensional range between $2$ and $9$.
    \item The dimensions of the centroid of associative  algebras of $5$-dimensional range between $2$ and $8$.
    \item The dimensions of the inner-derivation of associative  algebras of $5$-dimensional range between $2$ and $4$.
\end{itemize}
\end{corollary}
\begin{remark}
$\textbf{Dim(Inn(A))}\leq\textbf{Dim(Cent(A))}\leq\textbf{Dim(Der(A))}.$
\end{remark}

\end{document}